\pdfoutput=1
\documentclass[10pt]{article}
\usepackage{amsfonts,amsthm,amsmath,amssymb,url}

\usepackage{graphicx}

\newcommand{\reffig}[1]{Figure \ref{fig.#1}}

\makeatletter
\newcommand{\Pr@}{\operatorname{Pr}}
\newcommand{\E@}{\operatorname{E}}
\newcommand{\Var@}{\operatorname{Var}}
\renewcommand{\Pr}[1]{\ensuremath{\Pr@\left[{#1}\right]}}
\newcommand{\E}[1]{\ensuremath{\E@\left[{#1}\right]}}
\newcommand{\Var}[1]{\ensuremath{\Var@\left[{#1}\right]}}
\makeatother

\newcommand{\card}[1]{\ensuremath{\left\vert #1 \right\vert}}

\usepackage{url}

\newtheorem{theorem}{Theorem}
\newtheorem{corollary}[theorem]{Corollary}
\newtheorem{lemma}[theorem]{Lemma}

\newcommand{\reflem}[1]{Lemma \ref{lem.#1}}
\newcommand{\refthm}[1]{Theorem \ref{thm.#1}}

\newcommand{\refcor}[1]{Corollary \ref{cor.#1}}

\newcommand{\labellem}[1]{\label{lem.#1}}
\newcommand{\labelthm}[1]{\label{thm.#1}}

\newcommand{\labeleq}[1]{\label{eq.#1}}
\newcommand{\labelcor}[1]{\label{cor.#1}}

\def\complaint#1{}
\def\withcomplaints{
\newcounter{mycomplaints}
\def\complaint##1{\refstepcounter{mycomplaints}%
\ifhmode%
\unskip%
{\dimen1=\baselineskip \divide\dimen1 by 2 %
\raise\dimen1\llap{\tiny -\themycomplaints-}}\fi%
\marginpar{\tiny [\themycomplaints]: ##1}}%
}

\def\readRCS$#1: #2 #3 #4 #5${%
 \def\RCSfile{#2}%
 \def\RCSversion{#3}%
 \def\RCSdate{#4}%
}

\usepackage[colorlinks=false]{hyperref}
\usepackage{subfigure}
\withcomplaints

\author{Ruth Haas\thanks{Mathematics Department,
Smith College, Northampton, MA 01063. Email: {\tt rhaas@math.smith.edu}.} \and
Audrey Lee\thanks{
Department of Computer Science,
University of Massachusetts, Amherst, MA 01003. Email:
{\tt alee@cs.umass.edu}. Supported by an NSF graduate fellowship and NSF grant CCR-0310661.} \and
Ileana Streinu\thanks{
Computer Science Department, Smith College, Northampton, MA 01063.
Email: {\tt  streinu@cs.smith.edu}. Supported by  NSF grant CCR-0310661.} \and
Louis Theran\thanks{Department of Computer Science,
University of Massachusetts, Amherst, MA 01003. Email:
{\tt theran@cs.umass.edu}. Supported by  NSF grant CCR-0310661.}
}
\title{Characterizing Sparse Graphs by Map Decompositions}

\begin{document}
\maketitle

\begin{abstract}
A  {\bf map} is a graph that
admits an orientation of its edges so that each vertex has out-degree exactly
1. We characterize graphs which admit a decomposition
into $k$ edge-disjoint maps after: (1) the addition of {\it any} $\ell$
edges; (2) the addition of {\it some} $\ell$ edges.  These graphs are identified
with classes of {\it sparse} graphs; the results are also given in matroidal terms.

\end{abstract}

\section{Introduction and related work}
Let $G=(V,E)$ be a graph with $n$ vertices and $m$ edges.  In this paper, graphs
are multigraphs, possibly containing loops.  For a subset $V'\subset V$, we
use the notation $E(V')$ to denote the edges spanned by $V'$; similarly, $V(E')$
denotes the vertex set  spanned by $E'$.

A graph $G=(V,E)$ is {\bf $(k,\ell)$-sparse}, or simply sparse,\footnote{For brevity, we
omit the parameters $k$ and $\ell$ when the context is clear.}
 if no subset $V'$ of $n'$ vertices spans more than $kn'-\ell$ edges;
 when $m=kn-\ell$, we call the graph {\bf tight}.

Our interest in this problem stems from our prior work on pebble
game algorithms
\cite{streinu:lee:pebbleGames:EuroComb:2005,streinu:lee:theran:2005}.
The $(k,\ell)$-pebble game takes as its input a graph, and outputs
{\bf tight}, {\bf sparse} or {\bf failure} and an orientation of a
sparse subgraph of the input.  We had previously considered the
problem in terms of tree decompositions, suggesting the natural
range of $k\le \ell\le 2k-1$. In, fact, the pebble game generalizes
to the range $0\le\ell\le 2k-1$.  In this paper we examine the
graphs that the general pebble game characterizes.

 A  {\bf map} is a graph that
admits an orientation of its edges so that each vertex has
out-degree exactly 1.  This terminology and definition is due to
Lov{\'{a}}sz \cite{lovasz:combinatorial-problems:1979}. This class
of graphs is also known as the bases of the {\bf bicycle matroid}
\cite{oxley:matroidTheory:1992} or {\bf spanning pseudoforests}
\cite{gabow:westermann:matroidSums:1988}, where the equivalent
definition of having at most one cycle per connected component is
used.

Our choice of the former definition is motivated by the pebble game
algorithms. In the $(k,0)$-pebble game, the output orientation of a
tight graph has out-degree exactly $k$ for every vertex. The
motivation for studying the pebble game was to have a good algorithm
for recognizing sparse and tight graphs. These compute an
orientation of a sparse graph that obeys a specific set of
restrictions on the out degree of each vertex.

The focus of this paper is the class of graphs that decompose into
$k$ edge-disjoint {\bf maps} after the addition of $\ell$ edges;  we
call such a graph a {\bf $k$-map}.  Our goal is to extend the
results on adding $\ell-k$ edges to obtain $k$ edge-disjoint
spanning trees \cite{haas:arboricityGraphs:2002} to the range $0\le
\ell\le k-1$. A theorem of
\cite{streinu:lee:pebbleGames:EuroComb:2005} identifies the graphs
recognized by the $(k,\ell)$-pebble game as $(k,\ell)$-sparse
graphs.

\begin{figure}[htbp]
\centering 
\subfigure[]{\label{fig.k4}\includegraphics[height=.75 in]{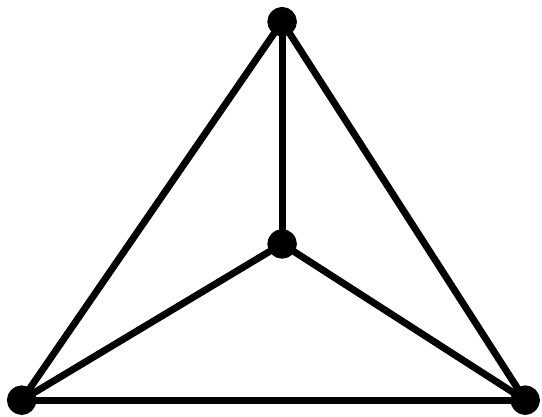}}
\hspace{.3in}
\subfigure[]{\label{fig.2-map}\includegraphics[height=.75 in]{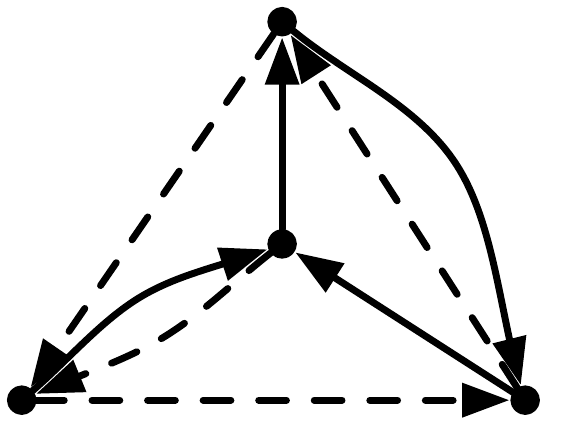}}
\hspace{.3in}
\subfigure[]{\label{fig.2-map-b}\includegraphics[height=.75 in]{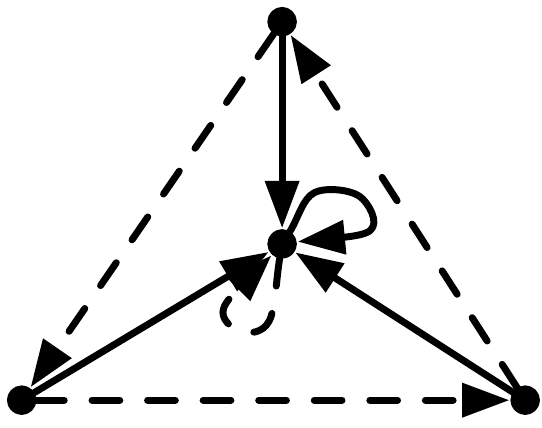}}
\caption{Adding any two edges to $K_4$ results in two maps.  }
\label{fig.add-to-k4}
\end{figure}

The complete graph
$K_4$  in \reffig{k4} is $(2,2)$-tight; i.e., adding
any two edges to $K_4$ we obtain a $2$-map.
The graphs in \reffig{2-map} and \reffig{2-map-b}
are obtained by adding two edges to $K_4$; the edges are
dashed and oriented to show a decomposition into two maps.

White and Whiteley \cite{whiteley:Matroids:1996} observe the
matroidal properties of sparse graphs for $0\le \ell\le 2k-1$  in
the context of bar-and-joint rigidity for frameworks embedded on
surfaces \cite{whiteley:unionMatroids:1988}. In
\cite{szego:constructive-sparse:egres-2003}, Szeg{\H{o}}
characterized exactly when tight graphs exist.

We also state our results in the context of matroid truncations. If
$\mathcal{M}=(E,\mathcal{I})$ is a matroid given by its independent
sets, then the {\bf truncation} of $\mathcal{M}$  is the matroid
$(E,\{E'\in \mathcal{I} : \card{E'}\le k\})$, for some nonnegative
integer $k$.  See, e.g., \cite{oxley:matroidTheory:1992} for a more
complete treatment of the topic of truncations.

The connection between sparse graphs and decompositions into
edge-disjoint spanning trees has been extensively studied.  The
classical results of Tutte
\cite{tutte:decomposing-graph-in-factors:1961} and Nash-Williams
\cite{nash-williams:edge-disjoint-spanning-trees:1961} show the
equivalence of $(k,k)$-tight graphs and graphs that can be
decomposed into $k$ edge-disjoint spanning trees; such a graph is
called a {\bf $k$-arborescence}.   A theorem of Tay
\cite{tay:thesis:1980,tay:rigidityMultigraphs-I:1984} relates such
graphs to generic rigidity of bar-and-body structures in arbitrary
dimension.

The particular case in which $k=2$ and $\ell=3$ has an important
application in rigidity theory: the minimally $(2,3)$-sparse graphs,
known as Laman graphs, correspond to minimally generically rigid
bar-and-joint frameworks in the plane \cite{laman:Rigidity:1970}.
Crapo \cite{crapo:genericRigidity:1996} showed the equivalence of
Laman graphs and those graphs that have a decomposition into $3$
edge-disjoint trees such that each vertex is incident to exactly $2$
of the trees; such a decomposition is called a {\bf $3T2$
decomposition}.

Of particular relevance to our work are results of Recski
\cite{recski:network-I:1984,recski:network-II:1984} and Lovasz and
Yemini \cite{lovasz:yemini:genericRigidity:1982}, which identify
Laman graphs as those that decompose into two spanning trees after
doubling any edge. In
\cite{hendrickson:thesis:1991,hendrickson:uniqueRealizability:1992}
Hendrickson characterized Laman graphs in terms of the existence of
certain bipartite matchings.  Stated in the terminology of this
paper, the results of \cite{hendrickson:thesis:1991} show the Laman
graphs are precisely those that decompose into $2$ edge-disjoint
maps after any edge is quadrupled.

The most general results linking sparse graphs to tree
decompositions are found in Haas \cite{haas:arboricityGraphs:2002},
who shows the equivalence of sparsity, adding $\ell-k$ edges to
obtain a $k$-arborescence, and $\ell Tk$ decompositions for the case
where $k\le \ell\le 2k-1$. Our results provide an analog of the
first equivalences in terms of graphs which decompose into $k$
edge-disjoint maps.

Another decomposition theorem involving sparse graphs is due to
Whiteley, who proved in \cite{whiteley:unionMatroids:1988} that for
the range $0\le \ell\le k-1$, the tight graphs are those that can be
decomposed into $\ell$ edge-disjoint spanning trees and $k-\ell$
edge-disjoint maps.

\section{Our Results}
Our results characterize the graphs which admit a decomposition into
$k$ edge-disjoint maps after adding $\ell$ edges.  Since the
focus of this paper is on the families of matroidal sparse graphs,
we assume that $0 \le \ell \le 2k-1$ unless otherwise stated.

First we consider the case in which we may add {\it any} $\ell$ edges, including
multiple edges and loops, to $G$.  Let $K_n^{k,2k}$ be the complete graph 
on $n$ vertices with $k$ loops on each vertex and edge multiplicity $2k$.  
It is easily seen that any sparse graph is a subgraph of $K_n^{k,2k}$, and 
we assume this in the following discussion.

\begin{theorem}
Let $G=(V,E)$ be a graph on $n$ vertices and $kn-\ell$ edges.
The following statements are equivalent:
\begin{enumerate}
\item $G$ is $(k,\ell)$-sparse (and therefore tight).
\item Adding {\it any} $\ell$ edges  from $K_n^{k,2k}-G$ to $G$ results in a $k$-map.
\end{enumerate}
\labelthm{maps-after-adding-any}
\end{theorem}

\refthm{maps-after-adding-any} directly generalizes the
characterization of Laman graphs in \cite{hendrickson:thesis:1991}.
It also generalizes the results of Haas
\cite{haas:arboricityGraphs:2002} to the range $0\le \ell\le k -1$.

As an application of \refthm{maps-after-adding-any} we obtain the
following decomposition result.

\begin{corollary}
Let $0\le \ell\le k$.  Let $G$ be a graph with $n$ vertices and $kn-\ell$ edges.
The following statements are equivalent:
\begin{enumerate}
\item $G$ is the union of
$\ell$ edge-disjoint spanning trees and $k-\ell$ edge-disjoint maps.
\item Adding {\it any} $\ell$ edges to $G$ results in a $k$-map.
\end{enumerate}
\labelcor{maps-and-trees}
\end{corollary}

We also characterize the graphs for which there are {\it some} $\ell$ edges that
can be added to create a $k$-map.

\begin{theorem}
Let $G=(V,E)$ be a graph on $n$ vertices and $kn-\ell$ edges.
The following statements are equivalent:
\begin{enumerate}
\item $G$ is $(k,0)$-sparse.
\item There is some set of $\ell$ edges, which when added to $G$
results in a $k$-map.
\end{enumerate}
\labelthm{add-some-edges}
\end{theorem}

Stating \refthm{add-some-edges} in matroid terms, we obtain the following.

\begin{corollary}
Let  $\mathcal{N}_{k,\ell}$ be the family of graphs $G$ such that
$m=kn-\ell$ and $G$ is $(k,0)$-sparse.  Then $\mathcal{N}_{k,\ell}$
is the class of bases of a matroid that is a truncation of the
$k$-fold union of the bicycle matroid. \labelcor{truncation}
\end{corollary}

Generalizing \refthm{maps-after-adding-any} and \refthm{add-some-edges} we have the
following theorem.

\begin{theorem}
Let $G=(V,E)$ be a graph on $n$ vertices and $kn-\ell-p$ edges and let $0\le\ell+p\le 2k-1$.
The following statements are equivalent:
\begin{enumerate}
\item $G$ is $(k,\ell)$-sparse.
\item There is some set $P$ of $p$ edges which when added to
$G$ results in a graph $G'=(V,E\cup P)$, such that adding any $\ell$
edges to $G'$ (but no more than $k$ loops per vertex) results in a
$k$-map.
\end{enumerate}
\labelthm{some-any}
\end{theorem}

In the next section, we provide the proofs.

\section{Proofs}
The proof of \refthm{maps-after-adding-any} relies on the following lemma.

\begin{lemma}\labellem{maps-are-sparse}
A graph $G$ is a $k$-map if and only if $G$ is $(k,0)$-tight.
\end{lemma}
\begin{proof}
Let $B_k(G)=(V_k,E,F)$ be the bipartite graph with one vertex class indexed by $E$ and
the other by $k$ copies of $V$.  The edges of $B_k(G)$ capture the incidence structure of
$G$.  That is, we define $F=\{ v_ie : e=vw, e\in E, i=1,2,\ldots,k\}$; i.e., each edge vertex in $B$
is connected to the $k$ copies of its endpoints in $B_k(G)$.  \reffig{bipartite} shows $K_3$
and $B_1(K_3)$.

\begin{figure}[htbp]
\centering 
\includegraphics[height=.75 in]{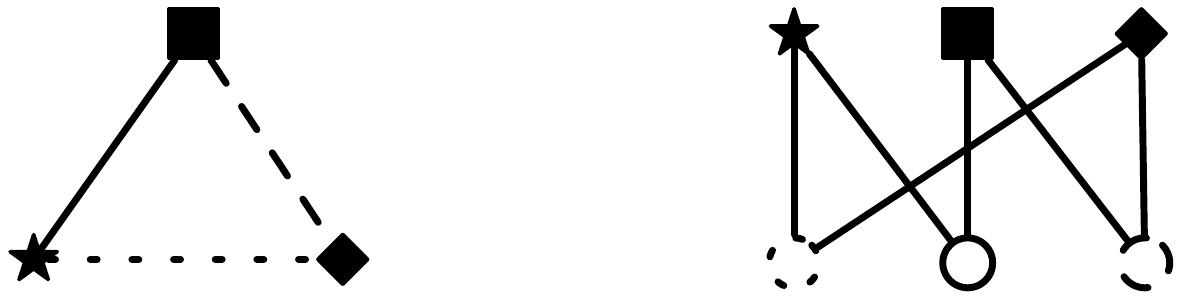}
\caption{$B_1(K_3)$
is shown on the right with the one copy of $V$ at the top. The style of line of the
edges on the left matches the style of line of the vertex in the bipartite graph corresponding
to that edge. }
\label{fig.bipartite}
\end{figure}

\begin{figure}[htbp]
\centering 
\includegraphics[height=.75 in]{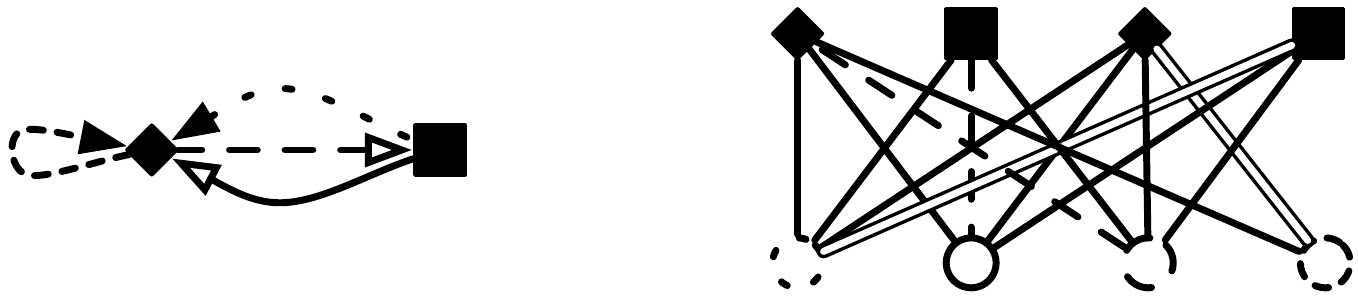}
\caption{$B_2(G)$ for the graph $G$ on the left
is shown on the right with the two copies of $V$ at the top.
$G$ is a 2-map; one possible decomposition is indicated
by the orientation of the edges and the style of arrow heads.  The matching
corresponding to this decomposition is indicated in the bipartite graph by dashed
and doubled edges.
}
\label{fig.bipartite-2-map}
\end{figure}

Observe that for $E'\subset E$, $N_{B_k(G)}(E')$, the neighbors of $E'$ in $B_k(G)$ of $E'$,  are
exactly the $k$ copies of the vertices of the subgraph spanned by $E'$ in $G$.  It follows
that
\begin{eqnarray}
\card{N_{B_k(G)}(E')}=k\card{V_G(E')}\ge \card{E'}
\labeleq{map-hall-condition}
\end{eqnarray}
holds for all $E'\subset E$ if and only  if $G$ is $(k,0)$-sparse. Applying Hall's theorem
shows that $G$ is $(k,0)$-tight if and only if $B_k(G)$ contains a perfect matching.

The edges matched to the $i$th copy of $V$ correspond to
the $i$th map in the $k$-map, as
shown for a 2-map in \reffig{bipartite-2-map}.
Orient each edge away from the vertex to which it is matched.  It follows that
each vertex has out degree one in the spanning subgraph matched to
each copy of $V$ as desired.
\end{proof}

\begin{proof}[Proof of \refthm{maps-after-adding-any}]
Suppose that $G$ is tight, and let $G'$ be the graph obtained by
adding any $\ell$ edges to $G$ from $K_n^{k,2k}-G$.
Then $G'$ has $kn$ edges; moreover $G'$ is $(k,0)$-sparse
since at most $\ell$ edges were added to the span of any subset $V'$
of $V$ of size at least $2$. Moreover,  since the added edges 
came from $K_n^{k,2k}$, they do not violate sparsity on single-vertex
subsets. It follows from \reflem{maps-are-sparse} that $G'$ can be
decomposed into $k$ edge-disjoint maps.

For the converse, suppose that $G$ is not tight.  Since $G$ has $kn-\ell$ edges, $G$ is
not sparse.  It follows that $G$ contains a subgraph $H=(V',E')$ such that
$\card{E'}\ge k\card{V'}-\ell+1$.  Add $\ell$ edges to the span of $V'$ to form $G'$.  By
construction $G'$ is not $(k,0)$-sparse; $V'$ spans at least $k\card{V'}+1$ edges in
$G'$.  Applying \reflem{maps-are-sparse} shows that $G'$ is not a $k$-map.
\end{proof}

\begin{proof}[Proof of \refcor{maps-and-trees}]
The equivalence of tight graphs for $0\le \ell\le k$ and the
existence of a decomposition into $\ell$ edge-disjoint spanning
trees and $(k-\ell)$ edge-disjoint maps is shown in
\cite{whiteley:unionMatroids:1988}. By
\refthm{maps-after-adding-any}, the tight graphs are exactly those
that decompose into $k$ edge-disjoint maps after adding any $\ell$
edges.
\end{proof}

\begin{proof}[Proof of \refthm{add-some-edges}]  By hypothesis, $G$ is $(k,0)$-sparse
but not tight.   By a structure theorem of
\cite{streinu:lee:pebbleGames:EuroComb:2005}, $G$ contains a single
maximal subgraph $H$ that is $(k,0)$-tight.  It follows that any
edge with at least one end in $V-V(H)$ may be added to $G$ without
violating sparsity. Adding $\ell$ edges inductively produces a tight
graph $G'$ as desired. Apply \reflem{maps-are-sparse} to complete
the proof.
\end{proof}

\begin{proof}[Proof of \refcor{truncation}]
Let $\mathcal{M}_k$ be the $k$-fold union of the bicycle matroid.  The
bases of $\mathcal{M}_k$ are exactly the $k$-maps.
Combining this with \refthm{add-some-edges} shows that $G\in \mathcal{N}_{k,\ell}$ if and
only if $G$ is independent in $\mathcal{M}_k$ and $\card{E(G)}=kn-\ell$ as desired.
\end{proof}

\begin{proof}[Proof of \refthm{some-any}]
Suppose that $G$ is sparse.  Since $G$ has $kn-\ell-p$ edges, $G$
does not contain a spanning $(k,\ell)$-tight subgraph.  Hence there
exist vertices $u$ and $v$ not both in the same $(k,\ell)$-tight
subgraph.  Add the edge $uv$.  Inductively add $p$ edges this way.
The resulting graph $G'$ is $(k,\ell)$-tight. By
\refthm{maps-after-adding-any}, adding any $\ell$ edges to $G'$
results in a $k$-map.

Now suppose that $G$ is not sparse.  As in \refthm{maps-after-adding-any},
there is no set of edges that can be added to $G$ to create a $(k,\ell)$-tight $G'$,
which proves the converse.
\end{proof}

\section{Conclusions and open problems}
\begin{figure}[htbp]
\centering 
\includegraphics[width=4 in]{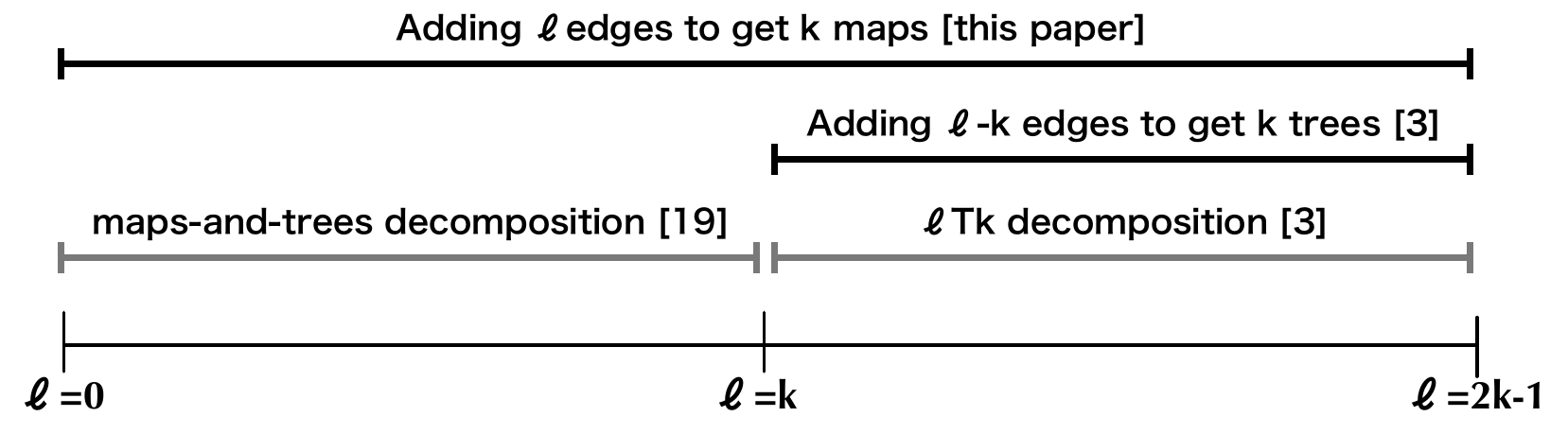}
\caption{Equivalent characterizations of sparse graphs in terms of decompositions and
adding edges.}
\label{fig.range-of-l-with-equivs}
\end{figure}

We characterize the graphs for which adding $\ell$ edges results in
a $k$-map.  These results are an analog to those of Haas
\cite{haas:arboricityGraphs:2002} using $k$-maps as the primary
object of study. In this setting, we obtain a uniform
characterization of the tight graphs for all the matroidal values of
$\ell$. \reffig{range-of-l-with-equivs} compares our results to
other characterizations of sparse graphs. In this paper we extend
the results of \cite{haas:arboricityGraphs:2002} to a larger range
of $\ell$.  While we do not have an analog of $\ell Tk$
decompositions for the new $0\le \ell\le k-1$ range, we do show the
equivalence of adding $\ell$ edges and the existence of a
decomposition into maps and trees.

In \cite{haas:arboricityGraphs:2002}, there are two additional types
of results: inductive sequences for the sparse graphs and the $\ell
T k$ decompositions.  Describing an analog of $\ell Tk$
decompositions for the maps-and-trees range of $\ell$ is an open
problem.

Lee and Streinu describe inductive sequences based on the pebble
game for all the sparse graphs in
\cite{streinu:lee:pebbleGames:EuroComb:2005}, but these do not give
the  explicit decomposition shown to exist in
\refcor{maps-and-trees}. Providing this decomposition explicitly
with an inductive sequence, as opposed to algorithmically as in
\cite{gabow:westermann:matroidSums:1988}, is another open problem.
The theorem of \cite{whiteley:unionMatroids:1988} used in the proof
of \refcor{maps-and-trees} is formulated in the setting of matroid
rank function and does not describe the decomposition.

\bibliography{maps}

\end{document}